\documentclass[11pt,reqno]{amsart}
\usepackage{amssymb,epsfig}

\subjclass{Primary: 37B20; Secondary: 37C45}

\newtheorem{proposition}{Proposition}[section]

\newtheorem{theorem}{Theorem}
\newtheorem{corollary}{Corollary}
\newenvironment{remark}{\noindent {\bf Remark:}}{}
\renewenvironment{proof}{\noindent {\bf Proof.}}{ \hfill\qed\\ }
\newenvironment{proofof}[1]{\noindent {\bf Proof of #1.}}{ \hfill\qed\\ }

\def\l{\lambda}

\def\ls{{\lambda^s}}
\def\lu{{\lambda^u}}
\def\Ls{{\Lambda^s}}
\def\Lu{{\Lambda^u}}
\def\A{\mathcal{A}}
\def\Z{\mathcal{Z}}

\def\N{\mathbb{N}}
\def\R{\mathbb{R}}
\def\T{\mathbb{T}}
\def\vrho{\rho}
\def\z{\zeta}
\def\hx{\hat{x}}
\def\hy{\hat{y}}
\def\O{\Omega}
\DeclareMathOperator{\supp}{supp}

\def\M{\mathcal{M}}
\let\al=\alpha
\let\eps=\varepsilon
\let\petit=\null

\def\e{\varepsilon}

\def\Z{\mathcal{Z}}

\def\1{{{\mathit 1} \!\!\>\!\! I} }
\renewcommand{\phi}{\varphi}
\newcommand{\br}{B(x,r)} 
\renewcommand{\liminf}{\mathop{{\underline {\hbox{{\rm lim}}}}}}
\renewcommand{\limsup}{\mathop{{\overline {\hbox{{\rm lim}}}}}}
\renewcommand{\eqref}[1]{{(\ref{#1})}}
\DeclareMathOperator{\esssup}{esssup}
\newcommand{\eqdef}{:=}

\begin{document}

\title{Recurrence and Lyapunov exponents}
\author{B.~Saussol \and S.~Troubetzkoy \and S.~Vaienti}
\address{LAMFA -- CNRS umr 6140, Universit\'e de Picardie Jules Verne, 
33, rue St Leu, F-80039 Amiens cedex 1, France}
\email{benoit.saussol@u-picardie.fr}
\urladdr{http://www.mathinfo.u-picardie.fr/saussol/}
\address{Centre de physique th\'eorique; 
Institut de math\'ematiques de Luminy; 
F\'ederation de Recherches des Unit\'es de Math\'ematique de Marseille
and Universit\'e de la M\'editerran\'ee, 
Luminy, Case 907, F-13288 Marseille Cedex 9, France.}
\email{troubetz@iml.univ-mrs.fr}
\urladdr{http://iml.univ-mrs.fr/{\lower.7ex\hbox{\~{}}}troubetz/}
\address{Phymat, Universit\'e de Toulon et du Var; Centre de physique th\'eorique and 
F\'ederation de Recherches des Unit\'es de
Math\'ematique de Marseille, CNRS Luminy, Case 907, F-13288 Marseille Cedex 9, France}
\email{vaienti@cpt.univ-mrs.fr}
\urladdr{http://www.cpt.univ-mrs.fr/}
\keywords{return time, Lyapunov exponents}
\begin{abstract}
We prove two inequalities between the Lyapunov exponents of a diffeomorphism
and its local recurrence properties.  We give examples showing that
each of the inequalities is optimal.
\end{abstract}
\maketitle

\bibliographystyle{plain}

\section{Introduction}
Given an ergodic map $f$ of a measure space $(M,\mathcal{B},\mu)$ which is also a metric
space and a measurable partition $\A$ of $M$
consider the ball $\br$ of radius
$r$ around the point $x$ and 
the $n$-cylinder around $A_n(x)\in {\A}^n := \bigvee_{j=0}^{n-1}f^{-j}{\A}$.
Define the first return of any set $A \in \mathcal{B}$  by
$$
\tau(A)  :=  \min \left \{k>0\colon
f^k(A)\cap A \neq\emptyset \right \}.
$$
This quantity is  also called the Poincar\'e
recurrence of the set $A$. It was previously
introduced in two contexts:
\begin{itemize}
\item {\it The statistics of return times into small neighborhoods}. 
It was shown in \cite{hsv}, that for certain weakly hyperbolic maps
$\tau(A_n(x))$ 
grows
linearly with $n$ when the metric entropy of
$\mu$ is positive.  For systems
with good mixing properties, this allows to ignore the
contribution of those points which come back
early, the other points giving an asymptotic
distribution of exponential type $e^{-t}$. 
\item \textit{The recurrence (or {\rm Afraimovich-Pesin}) dimension}. 
In  \cite{a,psv}, the quantity $\tau(\br)$ was used as a set function in the
Caratheodory covering of a given space, thus
providing a parametrized family of Borel measures
with transition point (dimension) located at some
dynamical characteristic of the space (usually
the topological entropy).
\end{itemize}
Especially motivated by the first item, we
proved a general result on {\it any}
measurable dynamical systems with positive
entropy $h_{\mu}$. In \cite{stv} we showed
\begin{quote}
{\em If the partition $\A$ is finite or
countable and the entropy
$h_{\mu}(f|\A)$ is strictly positive, then 
\begin{equation}\label{complexity0}
\liminf_{n\to\infty}\frac{\tau(A_n(x))}{n}\ge 1,
\quad \text{for $\mu$-almost every $x$}.
\end{equation}}
\end{quote}
Under the additional assumption that $f$ is invertible the two sided case 
\begin{equation}\label{complexity}
\liminf_{m,n\to\infty}\frac{\tau(A_{m,n}(x))}{m+n}\ge 1,
\quad \text{for $\mu$-almost every $x$}.
\end{equation}
was proven in \cite{acs}. 
Here $x\in A_m^n(x)\in \A_m^n:=\bigvee_{j=-m}^n f^{-j}\A$.
In both cases equality holds if we assume the
specification property.

In \cite{stv}, we proposed a new
technique to compute the Lyapunov exponent for a
large class of  weakly hyperbolic maps of the interval.  The method
expressed the exponent in terms of local recurrence properties.
This is one result in the spirit of what
was recently called ``the thermodynamics of return
times'': reformulating
the statistical properties of dynamical systems in
terms of returns of points or sets in the
neighborhood of (or into) themselves. 

Comparing cylinders with balls 
naturally relates $r^{-1}$ with the growth
(in the expanding case) of the derivative
$D_xf^n$ around the center
$x$ and has  two immediate consequences: it
allows to compare a ball of radius $r$ with a
cylinder of order $n := n(r)$ and to compare
$r$ with (some) Lyapunov exponents $\Lambda$ of
the map {\it via} the identification: $r^{-1}\sim
|D_xf^n|\sim e^{n \Lambda}$. In \cite{stv} this idea was
worked out rigorously for 
piecewise monotone maps of an interval with a finite
number of branches and with bounded derivative
of $p$-bounded variation with
an invariant measures with positive entropy.
For this class we approximated a cylinder from
the inside with a ball, used the
monotonicity of the set function $\tau{(\cdot)}$
and invoked 
\eqref{complexity0} to show that:
\begin{eqnarray}
\liminf_{r\rightarrow 0}\frac{\tau(\br)}{-\log
r}\ge \frac{1}{\lambda_{\mu}} \quad  \quad \mu-\hbox{a.e.}\label{formula}
\end{eqnarray}
where $\lambda_{\mu}$ is the positive Lyapunov
exponent of the measure $\mu$ of positive metric
entropy. When the map satisfies the specification
property the existence of the limit in (\ref{complexity0}) implies the
existence of the limit in (\ref{formula}), and its equality 
with the inverse of the Lyapunov exponent.

In this paper we generalize formula
(\ref{formula}) for multidimensional
transformations, providing examples which show that our inequalities
are sharp and other examples where 
strict inequalities hold. We will
work with $C^{1+\alpha}$ diffeomorphisms of a
compact manifold $M$ endowed with an invariant
ergodic measure $\mu$ of positive metric entropy
or invertible maps which are locally $C^{1+\alpha}$ diffeomorphisms 
with reasonable singularity sets such as in the monograph~\cite{ks}.

Let
$\lambda_1 \le \lambda_2 \le \dots \le \lambda_n$ be the
Lyapunov exponents of the measure $\mu$. In this very
general setting, we will prove that the left hand side in
(\ref{formula}) will be bounded from below by
$\frac{1}{\Lambda^{u}}-\frac{1}{\Lambda^{s}}$,
where $\Lambda^{u}$ (resp.~$\Lambda^{s}$) is the
largest $(=\lambda_n)$ (resp.~smallest
$(=\lambda_1)$) Lyapunov exponent. 
The technique of the proof will consists
in using (\ref{complexity}) together with a
comparison of balls and cylinders in local
Lyapunov charts.  In order to get an upper bound
for $\limsup \frac{\tau(\br)}{-\log r}$, 
one needs  the additional properties
that the measure is hyperbolic (all the Lyapunov
exponents are non-zero) and what we call a
``nonuniform specification property'' which will
insure the existence of periodic points,
of a well specified period, following closely the
orbit of any other point $x$ in a big set up to
a certain time $n$ (see the next section
 for the
precise definition). Using this assumption we
will compare balls not with cylinders, but with
Bowen (dynamical) balls, which will
contain periodic points and which
will  return  therefore into themselves with a
controlled time. In this case we will get an
upper bound of the form
$\frac{1}{\lambda^{u}}-\frac{1}{\lambda^{s}}$,
 where
$\lambda^{u}$ is the smallest positive Lyapunov
exponent and $\lambda^{s}$ is the largest negative
Lyapunov exponent. These bounds become
particularly easy, and more accessible, for
two-dimensional compact hyperbolic sets (or more
generally for two-dimensional invariant sets
enjoying the non-uniform specification property
quoted above): in these cases the limit exists and
gives the difference of the inverses of the
positive and negative Lyapunov exponents.

We point out that the product of this difference with the metric entropy
gives the Hausdorff dimension of the measure
$\mu$ (Young's formula). Furthermore if the map has a
constant Jacobian, the knowledge of one Lyapunov
exponent will immediately gives the other.
The numerical computation of the quantity
$\frac{\tau(\br)}{-\log r}$ seems very efficient
and reliable when applied to plane attractors and
conformal repellers; the detailed numerical analysis will be
published elsewhere.

In the last section we provide several examples.
First we provide examples where the limit of
$\frac{\tau(\br)}{-\log r}$ exists and equals the bounds
of Theorem \ref{thm1}  showing that each of the inequalities 
is optimal. Furthermore we provide an example where 
both inequalities are strict. In this example the limit equals
a mixture of the $\lambda_i$.

As a byproduct of these results on return times we obtain some information
on the spectrum of Poincar\'e recurrence of an invariant measure 
introduced in \cite{acs}. It turns out that the quantity 
$\lim_{r\to0}\frac{\tau(\br)}{-\log r}$ constitutes one of the component
of the pointwise dimension associated to the Caratheodory structure 
attached to the recurrence dimension. See Section~\ref{ap} for details.

\section{Basic definitions}

Throughout the article we suppose that
$M$ is an $n$--dimensional Riemannian manifold and $f: M \to M$ 
a $C^{1+\alpha}$ diffeomorphism
\footnote{Again, the results hold also for piecewise differentiable maps 
with reasonable singularity sets such as those considered in \cite{ks}.}.
A $f$--invariant measure $\mu$ is
called {\em hyperbolic} if all its Lyapunov exponents are non-null
and {\em partially hyperbolic} if at least one of its positive Lyapunov 
exponents and at least one of its negative Lyapunov exponents are non-null.
Suppose $\mu$ is ergodic and 
$\l_1 \le \l_2 \le \cdots \le \l_n$
are the corresponding Lyapunov exponents. 
The Ruelle inequality implies that if the entropy $h_\mu(f)$ is positive
then $\mu$ is hyperbolic or partially hyperbolic. 
Thus if $\mu$ is ergodic and of positive entropy then we can set
$\Ls := \l_1, \Lu := \l_n, \ls :=  \max\{\l_i: \l_i < 0 \}$ and $\lu := \min\{
\l_i: \l_i > 0\}$.

Next we define Lyapunov charts on the support of an ergodic measure $\mu$. 
Let $\R^j$ be the $j$-dimensional Euclidean space. Let
$s := \# \{\l_i: \l_i < 0\}$,  
$c := \# \{\l_i: \l_i = 0\}$ and
$u := \# \{\l_i: \l_i > 0\}$. 
For $x :=(x^u,x^c,x^s) \in \R^u \times \R^c \times \R^s$ let 
$|x| := |x^u|_u+|x^c|_c+|x^s|_s$ where 
$|\cdot|_u$, $|\cdot|_c$ and $|\cdot|_s$ are the Euclidean norms on $\R^u$,
$\R^u$ and $\R^s$. Let $|x|_u := |x^u|_u$, $|x|_c := |x^c|_c$
and $|x|_s := |x^s|_s$.
We define a distance function
$\vrho(z,z') := |z -z'|$.
The closed disk of radius $\rho$ centered at $x$ in $\R^u$ is denoted by 
$R^u(x,\rho)$ and 
$R(\rho) := R^u(0,\rho) \times R^c(0,\rho) \times R^s(0,\rho).$
Fix $\eta > 0$. For $\mu$--almost every $x \in M$ there is a Borel function 
$q(x) := q(x,\eta) : M \rightarrow [1,\infty)$ 
and an embedding $\Phi_x(R(q^{-1}(x))) \rightarrow M$ such that the following 
hold:

{\ \ \it i)}
$q(f^{\pm 1} x) \le q(x) e^{\eta}$, i.e. 
$q(x)$ is $\eta$-slowly varying,

{\ \ \it ii)}
$\Phi_x0 = x$, $D\Phi_x(0)$ takes $\R^u$, $\R^c$ and $\R^s$ 
to  $E^u(x)$, $E^c(x)$ and $E^s(x)$ respectively,

{\ \ \it iii)}
let  $f_x := \Phi^{-1}_{fx} \circ f \circ \Phi_x$ then
\[
\begin{array}{rcll}
e^{\Ls - \eta}|v| \le &|Df_x(0)v|& \le e^{\ls + \eta}|v| 
&\hbox{ for } v \in \R^s, \\ 
e^{\lu - \eta}|v| \le &|Df_x(0)v|& \le e^{\Lu + \eta}|v| 
&\hbox{ for } v \in \R^u,\\
e^{-\eta}|v| \le &|Df_x(0)v|& \le e^{\eta}|v|
&\hbox{ for } v \in \R^c,
\end{array}
\]

{\ \ \it iv)}
$\hbox{Lip}(f_x-Df_x(0)) < \eta$, $\hbox{Lip}(f^{-1}_x-Df^{-1}_x(0)) < 
\eta$, in $R(q^{-1}(x))$,

{\ \ \it v)} 
$C^{-1} d (\Phi_x z,\Phi_x z^{\prime}) \le 
\vrho(z,z') \le q(x) d(\Phi_x z,\Phi_x z^{\prime})$ 
for  $z,z^{\prime} \in R(q^{-1}(x))$
for some universal positive constant C.

Next we define a nonuniform hyperbolic version of the specification
property which we call {\em nonuniform specification} for $(f,\mu)$.
For $\mu$-almost every $x$, any integer $m$, $n$, and any $\e > 0$ 
there exists  $K := K(\eta,\eps,x,m,n)$ such that:

i) the \emph{nonuniform} Bowen ball 
\[
\widetilde B_m^n(x,\eps) \eqdef \bigcap_{k=-m}^n f^{-k}B(f^kx,\eps q(f^kx,\eta)^{-2})
\]
contains a periodic point with period $p\le n+m+K$;

ii) the dependence of  $K$  on $m,n$ satisfies
\begin{equation}\label{constant!}
\lim_{\eta\to0}\limsup_{m,n \to \infty} K(\eta,\e,x,m,n)/(m+n) = 0.
\end{equation}

We remark that even if $\supp(\mu)$ is a compact (uniformly) 
hyperbolic set for $f$, then in general $q(\cdot)$ cannot be 
chosen constant, so that there is no simple relation between
Bowen balls and nonuniform Bowen balls,
except the obvious one 
\[
\widetilde B_m^n(x,\eps) \subset B_m^n(x,\eps).
\] 
However we have the 
following relation.
\begin{proposition}\label{pro1}
Assume that $\supp(\mu)$ is a hyperbolic set for $f$ and that
$f|_{\supp(\mu)}$ satisfies the usual definition 
of specification (see for example \cite{kh}). 
Then $(f,\mu)$ is non-uniformly specified as well.
\end{proposition}
Recall that by the Spectral Decomposition Theorem 
$\supp\mu$ can always be decomposed into disjoint 
closed sets, where for some iterate of the map the 
specification property holds.

\begin{proof}
Denote the Bowen balls by
$$
B_m^n(z,\eps) := \bigcap_{k=-m}^n f^{-k}B(f^kz,\eps).
$$
Remark that by uniform hyperbolicity, if $\eps>0$ is sufficiently small
there exists $a\in(0,1)$ such that for all $z\in\supp\mu$
and integer $q$ we have
\begin{equation}\label{binb}
B_{q}^{q}(z,\eps) \subset B(z,a^q).
\end{equation}

Let $x$ be such that for any $\eta>0$ we have $q(x,\eta)<\infty$. 
Note that this concerns $\mu$-almost every $x\in M$.

Let $\nu>0$. Fix $\eta>0$ so small that $a^\nu e^{2\eta}<1$.
We write for simplicity $q(\cdot):=q(\cdot,\eta)$.
Assume that $m$ and $n$ are such that $m+n$ is so large that
$$
(a^\nu e^{2\eta})^{n+m} < a\eps q(x)^{-2}.
$$
Consider a ``nonuniform'' Bowen ball 
$$
\widetilde B_m^n(x,\eps) := \bigcap_{k=-m}^n f^{-k}B(f^kx,\eps q(f^kx)^{-2}).
$$
Our aim is to show that if $m$ and $n$ are sufficiently large then the 
nonuniform Bowen ball $\widetilde B_m^n(x,\eps)$ contains a periodic 
point with a sufficiently small period.
Observe that for any $k=-m,\ldots,n$ we have, by Equation (i) of
Lyapunov chart
$$
\eps q(f^kx)^{-2} 
\ge \eps q(x)^{-2} e^{-2\eta |k|} 
\ge \eps q(x)^{-2} e^{-2\eta (m+n)}.
$$
Let $q:=[\nu(m+n)]$.
By the specification property the Bowen ball $B_{m+q}^{n+q}(x,\eps)$
contains a periodic point, say $y$, with a period $p\le (m+n)(1+2\nu)+c$,
with the constant $c$ depending only on $\eps$.

For every $k=-m,\ldots,n$ we have $f^ky\in B_q^q(f^kx,\eps)$
hence by \eqref{binb}
$$
d(f^kx,f^ky) < a^q\le a^{\nu(n+m)-1} < \eps q(x)^{-2} e^{-2\eta(n+m)},
$$
the last inequality being satisfied since $m+n$ is sufficiently large.
Thus the periodic point $y$ belongs to $\widetilde B_m^n(x,\eps)$.
This shows that the function $K$ in \eqref{constant!} satisfies
$$
\limsup_{m+n \to \infty} K(\eta,\ell,\e,x,m,n)/(m+n) \le 2\nu.
$$ 
The proposition follows from the arbitrariness of $\nu$.
\end{proof}

\section{Statement of results}

\begin{theorem}\label{thm1}
If $\mu$ is an $f$--invariant, ergodic probability measure
with entropy $h_{\mu}(f) > 0$ then
$$ \frac{1}{\Lambda_\mu^u} - \frac{1}{\Lambda_\mu^s} \le 
\liminf_{r\to 0} \frac{\tau(B(x,r))}{\log 1/r} $$
for $\mu$--almost every $x$.

If $\mu$ is an $f$--invariant, ergodic probability measure
which is hyperbolic and $f|_{\supp(\mu)}$ satisfies the nonuniform 
specification property then
$$ \limsup_{r\to 0} \frac{\tau(B(x,r))}{-\log r} \le  
\frac{1}{\lambda_\mu^u} - \frac{1}{\lambda_\mu^s}$$
for $\mu$--almost every $x$.
\end{theorem}
The following corollaries show that inequalities in the theorem may
be optimal in some situation.
If $M$ is two dimensional
we have $\Lambda_\mu^u = \lambda_\mu^u$ and $\Lambda_\mu^s = \lambda_\mu^s$.
Furthermore whenever $h_{\mu}(f) > 0$ the Ruelle inequality
implies that $\mu$ is hyperbolic thus
the following corollary follows

\begin{corollary}\label{cor1}
If $M$ is two dimensional, $h_\mu(f) > 0$
and $f|_{\supp(\mu)}$ satisfies the nonuniform specification property then
$$ \lim_{r\to 0} \frac{\tau(B(x,r))}{-\log r} = 
\frac{1}{\lambda_\mu^u} - \frac{1}{\lambda_\mu^s} = \frac{1}{\Lambda_\mu^u} - \frac{1}{\Lambda_\mu^s}
$$
for $\mu$--almost every $x$.
\end{corollary}

\begin{corollary}\label{cor2}
If $f$ is a diffeomorphism and $\supp(\mu)$ is a compact 
locally maximal hyperbolic set for $f$ then 
$$ \limsup_{r\to 0} \frac{\tau(B(x,r))}{-\log r} \le  
\frac{1}{\lambda_\mu^u} - \frac{1}{\lambda_\mu^s}$$
for $\mu$--almost every $x$.

In particular if $M$ is two dimensional and $h_\mu(f)>0$
the limit exists and
equals the right hand side of this inequality.
\end{corollary}

\begin{proofof}{Corollary \ref{cor2}}
Using the Spectral Decomposition Theorem (see for example \cite{kh})
$\supp(\mu)$ is decomposed into disjoint closed sets
$\O_1,\ldots,\O_m$ which are permuted by $f$. If $k_i$ denotes
the smallest integer such that $f^{k_i}\O_i = \O_i$ then 
$f^{k_i}|_{\O_i}$ is topologically mixing.
In this situation it is well known that topological mixing is equivalent to the 
specification property, hence by Proposition~\ref{pro1} this implies that 
the nonuniform specification property is satisfied. 
The proof follows from Theorem~\ref{thm1} if we notice that 
whenever $x\in\O_i$ we have
$k_i \tau(B(x,r), f)  \le
k_i \tau(B(x,r), f|_{\supp{(\mu)}}) 
= \tau( B(x,r), f^{k_i}|_{\O_i})$
and $k_i\lu(f)=\lu(f^{k_i})$ and  $k_i\ls(f)=\ls(f^{k_i})$. 
\end{proofof}

\begin{proofof}{Theorem \ref{thm1}}
The proof of the first statement compares the ball 
$B(x,r)$ to a cylinder set of an
appropriate partition.  We need a finite
partition $\Z := \{\zeta_1,\zeta_2,\dots,\zeta_m\}$ of positive entropy
which satisfies the following property:
\begin{equation}\label{eq:1}
\mu(\{x \in M:  d(x,M \backslash \zeta(x)) < r\}) < cr
\end{equation} 
for all $r > 0$ for some positive constant $c$. Here $\zeta(x)$ is the
element of $\Z$ containing $x$.
The existence of such a partition was shown in \cite{bs}.

Fix $x$ such that for any $\eta>0$, $q(x,\eta)<\infty$.
Let $y \in M \backslash \z_0^{n}(x)$ satisfying 
\[
d(x,y) \le q^{-2}(x) e^{-(\Lu+ 3\eta)n}.
\]
Let $\hx := \Phi^{-1}_x(x)$ and $\hy := \Phi^{-1}_x(y)$. We remark that
$\hx=0$ by definition of Lyapunov charts.
The points $\hx$ and $\hy$ are close enough that we can apply 
Equation v) of Lyapunov charts, it yields 
$$\vrho(\hx,\hy)  \le \frac{e^{-(\Lu+3\eta)n}}{q(x)}.$$
For brevity we set $f^{(k)}_x := f_{f^{k-1}x} \circ \cdots \circ f_x$.
For $0 \le k \le n-1$ we
apply Equations ii)-iv) of Lyapunov charts to show that
$$\vrho \left (f^{(k)}_x \hx, f^{(k)}_x \hy \right ) \le 
\frac{e^{-(\Lu+3\eta)n} e^{(\Lu+2\eta)k}}{q(x)}.$$
We remark that a simple induction shows that Equation iv) can be used,
this induction follows from the fact that the right hand side of the
above equation is majorized by
$e^{-k\eta}q^{-1}(x) \le q^{-1}(f^kx)$. 

Again Equation v) of Lyapunov charts is applicable, it implies that
$$d(f^kx,f^ky) \le C  \frac{e^{-(\Lu+3\eta)n} e^{(\Lu + 2\eta)k}}{q(x)}.$$
Since $q(x) \ge 1$ this yields
$$d(f^kx,f^ky) \le C e^{-\eta n}.$$
Since $y \not \in \z_0^n(x)$ this implies that for some 
$k \in \{0,\dots, n-1\}$ the distance
between $f^kx$ and $M \setminus \z(f^kx)$ is less than 
$C e^{-\eta n}$.
This yields
\begin{eqnarray*}
&& \mu  \left (\left\{x \in M:  d(x,M \backslash \z_0^n(x))  \le  q^{-2}(x)e^{-((\Lu+ 3\eta)n} \right\}\right )\\
&& \quad \le n \cdot \max_{0 \le k \le n-1} 
\mu \left (\left\{ x \in M: d(f^kx,M \backslash \z(f^kx)) \le C  e^{-\eta n}
\right\}\right )\\
&&\quad \le n  \cdot c \cdot C  e^{-\eta n}.
\end{eqnarray*}
Here we used Equation \eqref{eq:1} and the $f$--invariance of $\mu$
in the last inequality.

By the Borel--Cantelli lemma we conclude that 
\begin{equation}\label{eq:2}
B(x, q^{-2}(x)e^{-(\Lu+ 3\eta)n}) \subset \zeta_0^n(x)
\end{equation}
for $\mu$--almost every $x$ for sufficiently large $n := n(x)$.
By considering $f^{-1}$ we have a similar statement for the backwards
direction
\begin{equation}\label{eq:3}
B(x, q^{-2}(x)e^{(\Ls - 3\eta)m}) \subset \zeta_m^0(x).
\end{equation}

Fix $r > 0$. Let $m_r$ be the largest integer such that  
$q^{-2}(x)e^{(\Ls - 3\eta)m} > r$ and $n_r$ the largest integer such
that  $q^{-2}(x)e^{-(\Lu+ 3\eta)n} > r$.

Combining Equations \eqref{eq:2} and \eqref{eq:3} with the definitions
of $m_r$ and $n_r$ yields
\begin{equation}\label{eq:4}
B(x, r) \subset \zeta_{m_r}^{n_r}(x).
\end{equation}

Thus
$$ 
\left ( \frac{1}{\Lu + 3\eta} - \frac{1}{\Ls - 3\eta}   \right )^{-1}
\liminf_{r \to 0} \frac{\tau(B(x,r))}{-\log r} \ge
 \liminf_{m,n \to \infty} \frac{\tau(\zeta_m^n(x))}{m+n}.
$$
In \cite{acs} 
it was shown that the right hand side 
is greater or equal to one. 
Since $\eta > 0$ is arbitrary this concludes the proof of the first
statement of the theorem.

We turn to the proof of the second statement.
For this proof it is more convenient to use Bowen balls.
Let 
$$
\widetilde B_m^n(x,\e) := \{y \in M: d(f^kx,f^ky)  < \e q^{-2}(f^kx) \quad
\forall -m \le k \le n\}.
$$
We need to ``reverse'' the inclusion in Equation \eqref{eq:4} using
nonuniform Bowen balls instead of partition elements.

Fix $\e > 0$ sufficiently small.  
We choose $m:=m(r)$ and $n:=n(r)$ the smallest possible integers such
that
\begin{equation}\label{eq:5}
C q^{-1}(x) \e \max(e^{(\ls + 3\eta)m} , e^{-(\lu - 3\eta)n} ) \le \frac{r}{2}.
\end{equation}
Consider $y \in\widetilde B_{m}^{n}(x,\e)$.
As before let $\hx := \Phi^{-1}_x(x)$ and $\hy := \Phi^{-1}_x(y)$.
In coordinates we write $\hx := (\hx^u,\hx^s)$ and 
$\hy := (\hy^u,\hy^s)$.

For the definition of Bowen balls it is clear that Equations i)-v)
of Lyapunov charts can be applied along the orbit segment for $f^{-m}x$ to 
$f^nx$. From Equation v) of Lyapunov charts we have
$\vrho(f^{(k)}_x\hx,f^{(k)}_x\hy) \le q^{-1}(f^kx) \e$ for any $-m\le k\le n$.
Equations ii)--iv) of Lyapunov charts imply that 
\[
\left\{
\begin{split}
\vrho(\hx^u,\hy^u) &\le q^{-1}(f^nx) \e e^{(-\lu + 2\eta)n},\\
\vrho(\hx^s,\hy^s) &\le q^{-1}(f^{-m}x) \e e^{(\ls + 2\eta)m}.
\end{split}
\right.
\]
The definition of $\vrho$ along with these two inequalities yields 
$$\vrho(\hx,\hy) \le q^{-1}(f^{-m}x) \e e^{(\ls + 2\eta)m} +
q^{-1}(f^nx) \e e^{-(\lu - 2\eta)n}.$$

Equation v) of Lyapunov charts implies that
\begin{eqnarray*}
d(x,y) &\le& C \e (q^{-1}(f^{-m}x)  e^{(\ls + 2\eta)m} +
q^{-1}(f^nx)  e^{-(\lu - 2\eta)n})\\
& \le & C q^{-1}(x) \e (e^{(\ls + 3\eta)m} + e^{-(\lu - 3\eta)n} )\\
& \le & r
\end{eqnarray*}
where the middle inequality use Equation i) of Lyapunov charts.
Thus we conclude that 
\begin{equation}\label{eq:6}
\widetilde B_{m}^{n}(x,\e) \subset B(x,r).
\end{equation}
Equations \eqref{eq:5} and \eqref{eq:6} imply that
\begin{equation}\label{eq:7}
\left ( \frac{1}{\lu - 3 \eta} - \frac{1}{\ls + 3 \eta}   \right )^{-1}
\limsup_{r \to 0} \frac{\tau(B(x,r))}{-\log r} \le
 \limsup_{m,n \to \infty} \frac{\tau(\widetilde B_m^n(x,\e))}{m+n}.
\end{equation}
By the nonuniform specification property we have
$$
\limsup_{n,m\to\infty}
\frac{\tau(\widetilde B_m^n(x,\e))}{m+n}
\le
1+\limsup_{n,m\to\infty}\frac{K(\eta,\eps,x,n,m)}{m+n}
\le 1+\delta(\eta),
$$
for some function $\delta$ such that $\lim_{\eta\to0}\delta(\eta)=0$.
We then combine this inequality with \eqref{eq:7}, and the conclusion 
follows from the arbitrariness of $\eta$.
\end{proofof}

\begin{remark}
In the case of endomorphisms, after the simplification consisting in 
ignoring possible stable directions, the proof of Theorem~\ref{thm1} 
can be carried out essentially in the same way.
Therefore for ergodic measure $\mu$ with entropy $h_\mu(f)>0$ it holds
\begin{equation}\label{thm1endo}
\frac{1}{\Lambda_\mu^u}  \le 
\liminf_{r\to 0} \frac{\tau(B(x,r))}{-\log r}
\quad\text{for $\mu$--almost every $x$.}
\end{equation}
Additionally, with the appropriate change to the notion of specification, 
namely considering ``forward'' Bowen balls $\widetilde B_0^n$, we get that 
if $\mu$ is an ergodic measure with all exponents positive and 
$f|_{\supp(\mu)}$ satisfies the nonuniform specification property then
\begin{equation}\label{thm1endo*}
\limsup_{r\to 0} \frac{\tau(B(x,r))}{-\log r} \le \frac{1}{\lambda_\mu^u} 
\quad\text{for $\mu$--almost every $x$.}
\end{equation}
\end{remark}

\section{Examples}
These examples show that the inequalities presented in the previous 
section are sometimes attained, but that still they may be strict.

The first example gives the optimality of the upper bound 
without the upper and lower bounds coinciding.

\begin{proposition}\label{pro:1}
Suppose that the linear automorphism $A := A_1 \times A_2$ of 
$\mathbb{T}^4$ is a direct product of two linear hyperbolic 
automorphisms $A_1, A_2$ of $\mathbb{T}^2$.

If $\mu_i$, $i=1,2$ are $A_i$--invariant ergodic probability measures 
with positive entropy then the $A$--invariant ergodic probability 
measure $\mu := \mu_1 \times \mu_2$ satisfies
$$ \lim_{r\to 0} \frac{\tau(B(x,r))}{-\log r} = 
\frac{1}{\lu} - \frac{1}{\ls}$$
for $\mu$--almost every $x$.
\end{proposition}

\begin{proof}
Since the limit in question depends only on the equivalence class of
 metrics
we suppose that the metric $d$ on $\T^4$ is defined by $d := \max(d_1,d_2)$
where for $i=1,2$ $d_i$ is the usual  metric on $\T^2$.
Any notation without a subscript pertains to the map $A$, while
the corresponding notation with a subscript pertains to the maps
$A_i$, i=1,2. For example
$\lambda^u_i$ and $\lambda^s_i$ denote the Lyapunov exponents of $A_i$.
Clearly $\lu = \min_{i} \lambda^u_i$ and $\lambda^s = \max_i \ls_i$.
We remark that these equalities are realized by the same $i$
since $\lambda^u_i + \lambda^s_i = 0$.

Let $x := (x_1,x_2)$ and notice that 
\begin{equation}\label{eq:9}
\tau(B(x,r)) \ge \tau_i(B(x_i,r)).
\end{equation}
Applying Corollary \ref{cor2} to $A_i$ yields
$$ \lim_{r\to 0} \frac{\tau_i(B(x_i,r))}{-\log r} = 
\frac{1}{\lambda^u_i} - \frac{1}{\lambda^s_i}$$
for $\mu_i$--almost every $x_i$.
Combining this with Equation \eqref{eq:9} yields
$$ \liminf_{r\to 0} \frac{\tau(B(x,r))}{-\log r} \ge 
\frac{1}{\lu} - \frac{1}{\ls}$$
for $\mu$--almost every $x$.
The conclusion follows by applying Corollary \ref{cor2} to the map $A$.
\end{proof}

The second example yields the optimality of the lower bound 
without the upper and lower bounds coinciding.

\begin{proposition}\label{pro:2}
Suppose that the linear automorphism $A := A_1 \times A_2$ of 
$\mathbb{T}^4$ is a direct product of two linear hyperbolic 
automorphisms $A_1, A_2$ of $\mathbb{T}^2$.

There exist some $A$--invariant ergodic probability measures $\mu$
with positive entropy such that
$$ \lim_{r\to 0} \frac{\tau(B(x,r))}{-\log r} = 
\frac{1}{\Lu} - \frac{1}{\Ls}$$
for $\mu$--almost every $x$.
\end{proposition}

\begin{proof}
We choose the metric as in the proof of Proposition~\ref{pro:1}
and keep the same notation.
Without loss of generality we can suppose that 
$(\Lu,\Ls)=(\lambda^u_1,\lambda^s_1)$.
Let $\mu_1$ be an ergodic $A_1$--invariant measure with positive
 entropy, and let $\mu_2$ be a Dirac mass at some fixed point $p$ of $A_2$.
The direct product $\mu:=\mu_1\times\mu_2$ is an $A$--invariant ergodic
probability measure with positive entropy.
For $\mu$-almost every $x$ we have $x_2=p$ which clearly implies
that $\tau(B(x,r)) = \tau_1(B(x_1,r))$ for any $r>0$.
Applying Corollary \ref{cor2} to $A_1$ finishes the proof.
\end{proof}
Notice that in this example the map has a zero entropy factor.

The third example is a map for which both inequalities in the theorem
are strict, and furthermore more than one of the unstable Lyapunov
exponents play a role. To simplify the exposition the following example 
is given by an expanding map of the torus, and not a diffeomorphism. 
Nevertheless we believe that the mechanism described there is widespread
for automorphisms of the torus.

\begin{proposition}
There exists a linear expanding map of the torus $\T^2$ such that 
$$ 
\frac{1}{\Lu} <
\lim_{r\to 0} \frac{\tau(B(x,r))}{-\log r} 
= \frac{2}{\Lu+\lu} < \frac{1}{\lu}
$$
for Lebesgue-almost all $x$.
\end{proposition}

\begin{proof}
Consider the linear expanding map $f:\T^2\to\T^2$ given by
$$
f(x) = Ax, \quad\text{ where } A= 
 \left ( \begin{matrix}
6 & 3 \\
3 & 3 
\end{matrix} \right )
$$
We first prove that for any $x$,
\begin{equation}\label{eq:limsup}
 \limsup_{r\to0} \frac{\log\tau(B(x,r))}{-\log r} 
\le \frac{2}{\Lu+\lu}
\end{equation}

The eigenvalues of $A$ are 
$e^\lu=\frac{9 - 3\sqrt{5}}2$ and $e^\Lu=\frac{9 + 3\sqrt{5}}2$
and the corresponding eigenvectors are 
$$
v^u =\left ( \begin{matrix}
1\\
\frac{- \sqrt{5} - 1}{2} 
\end{matrix} \right ) \hbox{ and }
V^u=\left ( \begin{matrix}
1\\
\frac{\sqrt{5} - 1}{2} 
\end{matrix} \right ).
$$
Let $m=\|v^u\|+\|V^u\|$. Fix $r>0$ and set 
$$
L(r)=[-r,r]v^u+[-r,r]V^u\subset \T^2.
$$
We have $L(r) \subset B(0,mr)$.
We want to find some integer $n$ such that
\begin{equation}\label{eq:11}
f^n L(r) \supset \T^2.
\end{equation}
Notice that the set
$$
A^n L(r) = 
[-e^{\lu n}r,e^{\lu n}r]v^u + [-e^{\Lu n}r,e^{\Lu n}r]V^u
$$
consists of a long and thin strip of width $2e^{\lu n}r$ parallel 
to the direction $V^u$, which is wrapped around the torus $2e^{\Lu n}r$ 
many times.
Let $\theta=\frac{\sqrt{5} - 1}{2}$, $c=1/\sqrt{1+\theta^2}$ and
let $p_i/q_i$ be the convergents of the continued fraction approximation
of $\theta$.

Let $\displaystyle n=\left\lceil \frac{-2\log r+\log{4c/3}}{\Lu+\lu} 
\right\rceil$.

We remark that $\theta$ has bounded quotients, 
more precisely $1<q_i/q_{i-1}\le 2$, thus there exists $i$ such that
$$
2ce^{-n\lu}/r\ge q_{i-1}>ce^{-n\lu}/r.
$$
For any $s\in \T^1$ there exists some integer $k=0,\ldots,q_i$ such that
\begin{equation}\label{eq:13}
\|s-k\theta\| < 1/q_{i-1} < ce^{n\lu}r.
\end{equation}
Let $R_\theta$ be the rotation by angle $\theta$ on $\T^1$.
Denote by 
$$
L^{uu}(r)=[-r,r]V^u\subset \T^2.
$$
Since $V^u=(1,\theta)$ the following inclusion holds (see Figure~\ref{fig:2})
$$
\{ (0,R_\theta^k(0))\colon k\in\N, |k|\le e^{n\Lu}r \}
\subset A^n (L^{uu}(r)).
$$
\begin{figure}
\begin{center}
\resizebox{8cm}{4cm}{\includegraphics{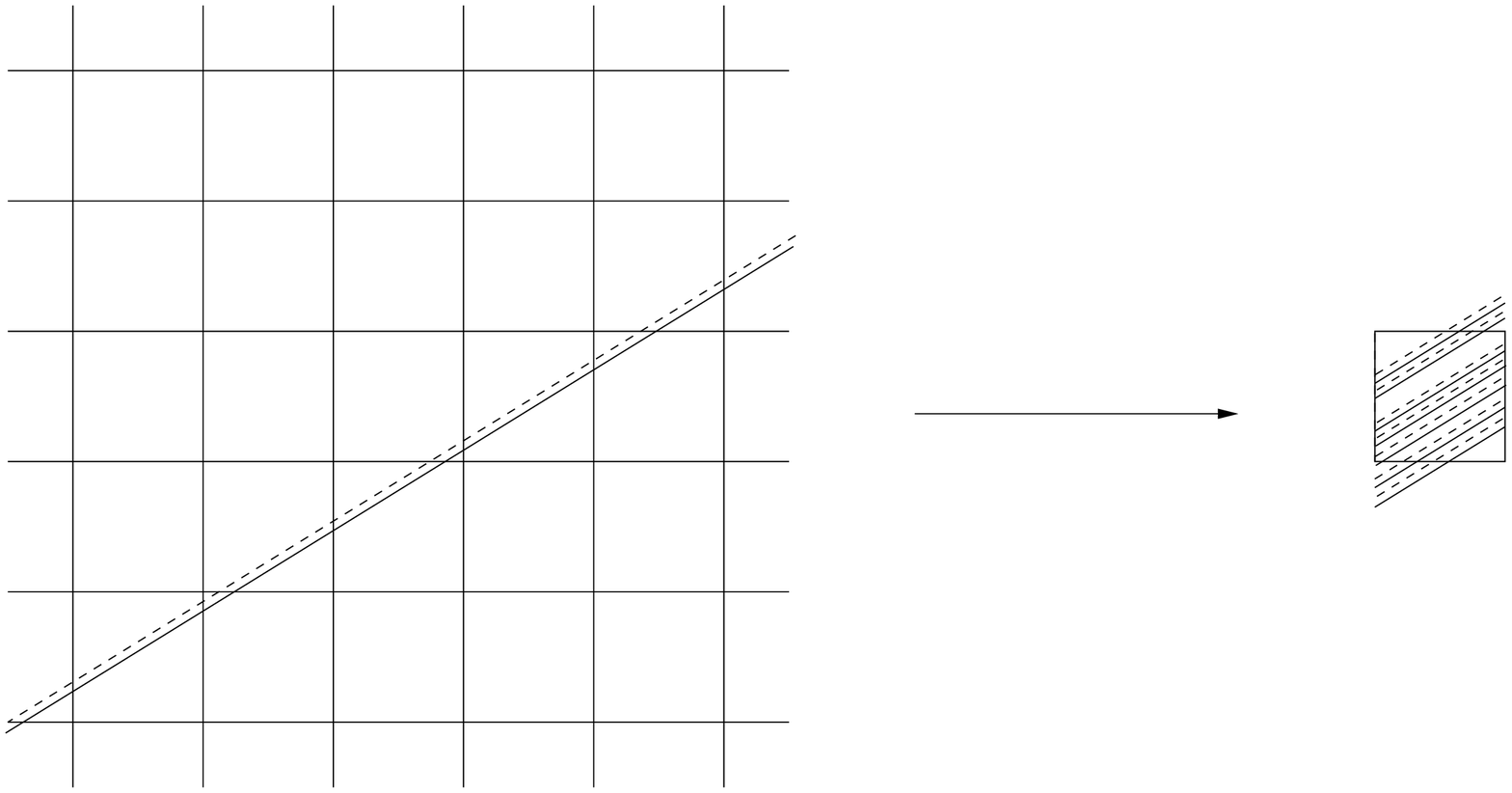}}
\put (-100,80){$A^nL^{uu}(r)$}
\put (-180,-10){$0$}
\put (-100,40){$\R^2$}
\put (-40,40){$\T^2$}
\end{center}
\caption{Intersection of the piece of strong unstable manifold 
$A^nL^{uu}(r)$ with $\{0\}\times \T^1$.}\label{fig:2}
\end{figure}

By our choice of $n$ we have
$q_i\le 2 q_{i-1}\le 4ce^{-n\lu}/r \le \lfloor e^{\Lu n}r\rfloor-1$,
hence by Equation~\eqref{eq:13} the set 
$P=\{R_\theta^k(0)\colon k=0,\ldots,\lfloor e^{\Lu n}r\rfloor-1\}$ 
is $ce^{n\lu}r$-dense in $\T^1$.
Therefore Equation \eqref{eq:11} holds (see Figure~\ref{fig:3}).
\begin{figure}
\begin{center}
\resizebox{2cm}{2cm}{\includegraphics{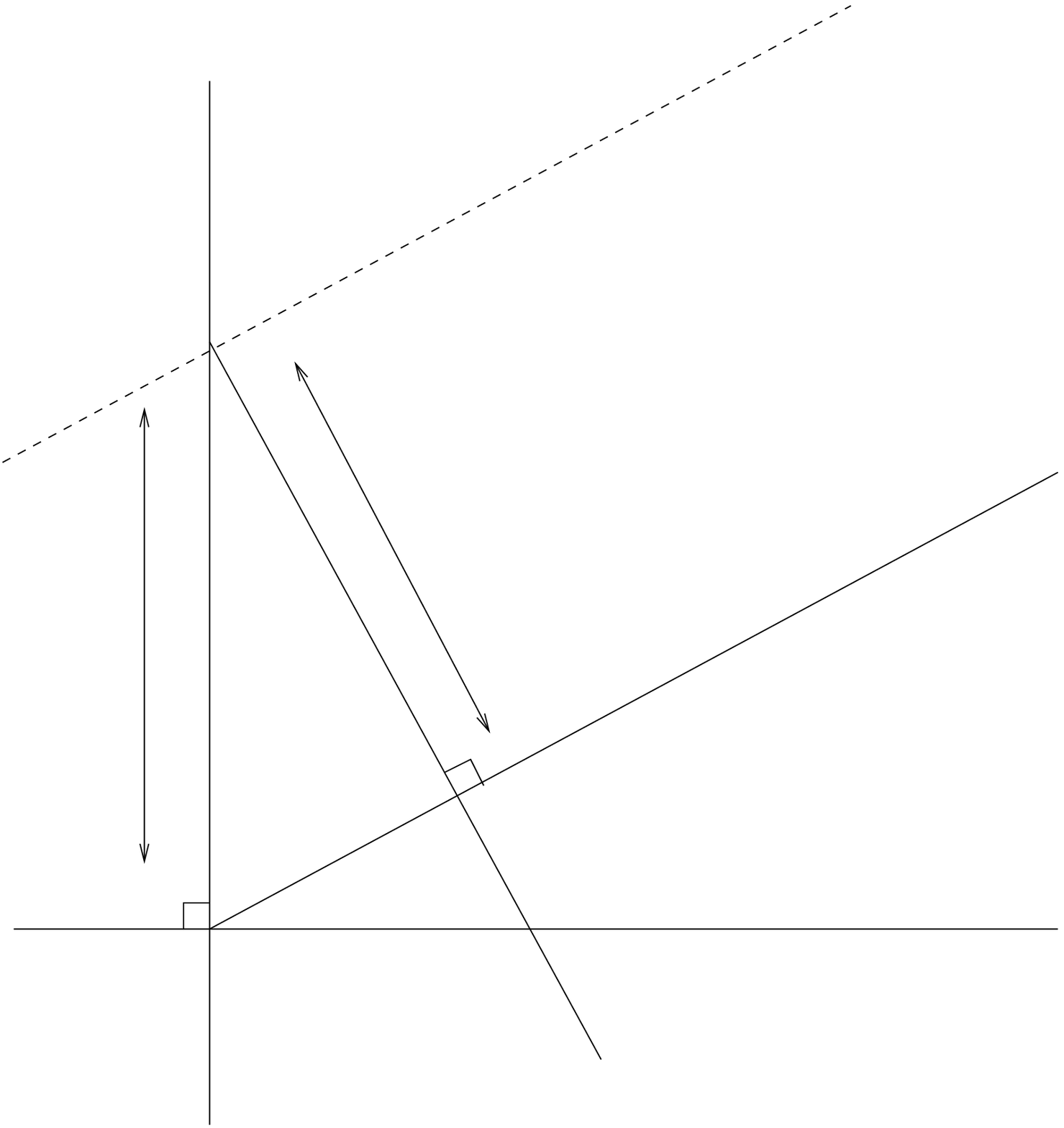}}
\put (-33,30){$e^{\lu n}r$}
\put (-82,20){$ce^{\lu n}r$}
\end{center}
\caption{$ce^{\lu n}r$-density in $\T^1$ of the set $P$ imply 
$e^{\lu n}r$-density of $A^nL^{uu}(r)$ in $\T^2$.}\label{fig:3}
\end{figure}
This certainly implies that $\tau(B(x,m r))\le n$ for any $x\in \T^2$,
and proves the upper bound \eqref{eq:limsup} for any points $x\in\T^2$.

We turn now to the almost sure lower bound.
Let $c=1-\| A^{-1}\|>0$.
Since the map $f$ is locally expanding the Closing Lemma implies that
for any $z\in\T^2$, real $r>0$ and integer $p$, if $d(z,f^pz)<r$ 
then there exists a periodic point $y$ with period $p$ such that 
$d(y,z)<r/c$.

Observe that if $p=\tau(B(x,r))$ then there exists $z\in \T^2$ 
such that $d(z,x)<r$ and $d(f^nz,x)<r$, thus in particular
$d(z,f^nz)<2r$, so that by the previous paragraph there exists 
a periodic point $y$ with period $p$ such that $d(y,z)<2r/c$. 
This implies that 
\[
d(x,F_p) \le (1+2/c)r, 
\quad\text{where } F_p=\{y\in\T^2\colon y=f^py\}.
\]
Let $a<(\det A)^{-1/2}$. Since the cardinality of $F_p$ equals
$\det (A^p-I)= (e^{p\lu}-1)(e^{p\Lu}-1)\le \det A^p$ we have
\[
\begin{split}
m(x\in\T \colon \tau(B(x,a^n))\le n) 
&\le \sum_{p=1}^n \sum_{y\in F_p} m( B(y,(1+2/c)a^n) )\\
&\le \sum_{p=1}^n (\det A)^p (1+2/c)^2a^{2n}\\
&\le (1+2/c)^2\frac{\det A}{\det A-1}(a^2\det A)^n.
\end{split}
\]
This is summable, hence by the Borel-Cantelli Lemma 
for Lebesgue almost all $x$
\[
\liminf_{n\to\infty} \frac{\tau (B(x,a^n))}{n} \ge 1.
\]
Taking $a$ arbitrarily close to 
\[
(\det A)^{-1/2}=
\exp \left(-\frac{\lambda^u+\Lambda^u}{2}\right)
\]
gives that 
\[
\liminf_{r\to 0} \frac{\tau (B(x,r))}{-\log r} 
\ge \frac{2}{\lambda^u+\Lambda^u}.
\]
\end{proof}
Notice that the second part of the proof may be followed in the more
general case of hyperbolic basic set of diffeomorphisms to gives
\[
\liminf_{r\to 0} \frac{\tau (B(x,r))}{-\log r} 
\ge \frac{\dim_H\mu}{h_{\textrm{top}}(f)}.
\]
This lower bound, although sometimes better than the one provided by
Theorem~\ref{thm1} (in dimension larger than $1+1$), is not optimal
even in dimension one. 
We conjecture that for any ergodic hyperbolic measure $\mu$ 
with nonzero entropy one has
\[
\liminf_{r\to 0} \frac{\tau (B(x,r))}{-\log r} 
\ge \frac{\dim_H\mu}{h_\mu(f)}.
\]
It can be shown that this bound is never worse than the one given 
by Theorem~\ref{thm1}.

\section{Application to the spectrum of recurrence dimension}\label{ap}
In this last section we propose to apply the main results to the 
recurrence dimension of measures introduced by Afraimovich {\it et al} 
in \cite{acs}.
We briefly recall the construction of the spectrum for Poincar\'e
recurrence of a measure.

Let $(M,f,\mu)$ be an ergodic measure preserving dynamical systems
on a compact manifold $M$.

For any $A\subset X$, $\al\in\R$ and $q\in\R$ we define
\begin{equation}\label{def-stat-sum*}
\M^{\petit \tau} (A,\al,q,\eps) \eqdef 
\inf
\sum_{i} \exp \big[-q \tau(B(x_i,\eps_i)) \big] \eps_{i}^{\al}\;,
\end{equation}
where the infimum is taken over all finite or 
countable collections of balls $B(x_i,\eps_i)$ such that 
$\bigcup_{i} B(x_{i},\eps_i)\supseteq A$ and $\eps_i\le\eps$.
 
The limit
$\M^{{\petit \tau}}(A,\al,q)\eqdef 
\lim_{\eps\rightarrow 0}\M^{{\petit \tau}}(A,\al,q,\eps)$
exists by monotonicity and we define
for any non-empty $A\subset X$ and any $q\in\R$,
\begin{equation}\label{pre-mesure}
\al_{{\petit \tau}}(A,q) \eqdef
\inf \{\al\colon\M^{{\petit \tau}}(A,\al,q)=0\}
= \sup\{\al\colon\M^{{\petit \tau}}(A,\al,q)=\infty\}
\end{equation}
with the convention that $\inf\emptyset =+\infty$ and $\sup\emptyset=-\infty$.
Notice that $\alpha_{\petit \tau}(A,0)$ is equal 
to the Hausdorff dimension of the set $A$.
Now we proceed to the definition of the spectra of measures, 
following~\cite{pesin}
\[
\al_{\petit{\tau}}^{\mu}(q)\eqdef \inf\{\al_{\petit{\tau}}(Y,q)\colon 
Y\text{ measurable }\subset X, \,\mu(Y)=1\}\;.
\]
This global quantity may be computed with the help of 
the corresponding pointwise dimension defined for any $x$ by
\[
d_{\mu,q}(x)= \liminf_{r\to0} \inf_{y\in B(x,r)}
 \frac{\log\mu(B(y,r))+q\tau(B(y,r))}{\log r}.
\]
A priori we cannot discard in general the infimum on $y\in B(x,r)$.
However, when $q=0$ we have the equality
\[
 d_{\mu,0}(x) = \liminf_{r\to0}\frac{\log\mu(B(x,r))}{\log r}.
\]
The well known equality $\dim_H\mu = \mu\text{-}\esssup d_{\mu,0}$
was generalized in~\cite{cs} for $q\neq0$:
\begin{equation}\label{eq:thm*}
 \al_{\petit{\tau}}^{\mu}(q) = \mu\text{-}\esssup d_{\mu,q}.
\end{equation}

\begin{corollary}\label{thm:2}
If $\mu$ is an $f$-invariant, ergodic measure with
 $h_\mu(f)>0$ then
\[
\begin{array}{cccl}
\al_{\petit{\tau}}^{\mu}(q) &\le& 
\dim_H\mu - q\left(\frac1{\Lambda_\mu^u}-\frac1{\Lambda_\mu^s}\right) 
& \text{ if $q\ge0$,}\\
\al_{\petit{\tau}}^{\mu}(q)&\ge&
\dim_H\mu - q\left(\frac1{\Lambda_\mu^u}-\frac1{\Lambda_\mu^s}\right)
& \text{ if $q\le0$.}
\end{array}
\]
If in addition $\mu$ is hyperbolic and $f|_{\supp\mu}$ 
satisfies the nonuniform specification property then 
$$\al_{\petit{\tau}}^{\mu}(q) \le \dim_H\mu -
q\left(\frac1{\lambda_\mu^u}-\frac1{\lambda_\mu^s}\right)
\text{ if $q\le0$.}
$$
\end{corollary}
The immediate corollary follows
\begin{corollary}\label{cor:3}
If $f$ is a surface diffeomorphism, $\mu$ is an $f$-invariant ergodic measure
with $h_\mu(f)>0$ and $\supp\mu$ is a compact locally maximal hyperbolic set 
for $f$ then for any $q\le 0$
\[
\begin{split}
\al_{\petit{\tau}}^{\mu}(q) &= 
\dim_H\mu - q\left(\frac1{\lambda_\mu^u}-\frac1{\lambda_\mu^s}\right)\\
&=
\left(1-\frac{q}{h_\mu(f)}\right)\dim_H\mu.
\end{split}
\]
\end{corollary}
\proofof{Corollary \ref{thm:2}}
Whenever $q\ge0$ we have
\[
\begin{split}
d_{\mu,q}(x) &\le 
\liminf_{r\to0} \frac{\log\mu(B(x,r))}{\log r}
-q\liminf_{r\to0} \frac{\tau(B(x,r))}{-\log r}\\
&\le
d_{\mu,0}(x) - q \left(\frac1{\Lambda_\mu^u}-\frac1{\Lambda_\mu^s}\right),
\end{split}
\]
for $\mu$-almost every $x$ by Theorem~\ref{thm1}.
Using Equation~\eqref{eq:thm*} applied to $d_{\mu,0}$ and 
$d_{\mu,q}$ yields the 
upper bound for $\al_{\petit{\tau}}^{\mu}(q)$ for $q\ge0$.

Suppose $q\le 0$. Notice that whenever 
$y\in B(x,r)$ we have $B(y,r)\subset B(x,2r)$, and thus
\[
\begin{split}
 d_{\mu,q}(x) 
&\ge \liminf_{r\to0} \frac{\log\mu(B(x,r))}{\log r}
-q\liminf_{r\to0} \frac{\tau(B(x,2r))}{-\log r} \\
&\ge
d_{\mu,0}(x) -q\left(\frac1{\Lambda_\mu^u}-\frac1{\Lambda_\mu^s}\right), 
\end{split}
\]
for $\mu$-almost every $x$ by Theorem~\ref{thm1}.
Again Equation~\eqref{eq:thm*} applied to $d_{\mu,0}$ and 
$d_{\mu,q}$ yields the lower bound for 
$\al_{\petit{\tau}}^{\mu}(q)$ when $q\le0$.

We prove the last part of the corollary.
If $q\le0$ we have
\[
\begin{split}
d_{\mu,q}(x)&\le \liminf_{r\to0} \frac{\log\mu(B(x,r))}{\log r}
-q \limsup_{r\to0}\frac{\tau(B(x,r))}{-\log r} \\
&\le d_{\mu,0}(x)-q\left(\frac1{\lambda_\mu^u}-\frac1{\lambda_\mu^s}\right), 
\end{split}
\]
for $\mu$-almost every $x$, by the second statement of
Theorem~\ref{thm1}.
Applying Equation~\eqref{eq:thm*} with $d_{\mu,0}$ and 
$d_{\mu,q}$ yields the result.
\endproof

\proofof{Corollary \ref{cor:3}}
Under these conditions, $f|_{\supp\mu}$ satisfies the nonuniform
specification property and the measure $\mu$ is hyperbolic. 
Since in addition 
$(\lambda_\mu^u,\lambda_\mu^s)=(\Lambda_\mu^u,\Lambda_\mu^s)$ 
the first equation follows from Corollary~\ref{thm:2}. 
The second equation is a consequence of Young's formula 
$$
\dim_H\mu = \frac{h_\mu(f)}{\lambda_\mu^u}- \frac{h_\mu(f)}{\lambda_\mu^s}
$$
for surface diffeomorphisms~\cite{y}.
\endproof


\end{document}